\documentclass[3p,10pt,times]{elsarticle}
\usepackage{amsthm,amssymb,amsmath}
\usepackage{hyperref}
\usepackage{subcaption}
\usepackage{mathtools}
\usepackage{xparse}
\usepackage{xcolor}
\usepackage[normalem]{ulem}

\DeclarePairedDelimiter{\ceil}{\lceil}{\rceil}

\theoremstyle{definition}
\newtheorem{theorem}{Theorem}
\newtheorem{example}[theorem]{Example}

\newtheorem*{innerdefinition}{Definition}
\newcounter{definition}

\makeatletter
\NewDocumentEnvironment{definition}{o}
 {%
  \IfValueTF{#1}
    {\innerdefinition[#1]\refstepcounter{definition}\def\@currentlabel{(#1)}}%
    {\innerdefinition}%
 }
 {%
  \endinnerdefinition
 }
\makeatother

\begin{document}
\baselineskip11pt

\begin{frontmatter}

\title{Extending discrete exterior calculus to a fractional derivative}

\author{Justin Crum}
\address{Program in Applied Mathematics, University of Arizona}
\ead{jcrum@math.arizona.edu}
\author{Joshua A. Levine}
\address{Department of Computer Science, University of Arizona}
\ead{josh@email.arizona.edu}
\author{Andrew Gillette}
\address{Department of Mathematics, University of Arizona}
\ead{agillette@math.arizona.edu}

\begin{abstract} 
Fractional partial differential equations (FDEs) are used to describe phenomena that involve a ``non-local'' or ``long-range'' interaction of some kind.
Accurate and practical numerical approximation of their solutions is challenging due to the dense matrices arising from standard discretization procedures.
In this paper, we begin to extend the well-established computational toolkit of Discrete Exterior Calculus (DEC) to the fractional setting, focusing on proper discretization of the fractional derivative.
We define a Caputo-like fractional discrete derivative, in terms of the standard discrete exterior derivative operator from DEC, weighted by a measure of distance between $p$-simplices in a simplicial complex.
We discuss key theoretical properties of the fractional discrete derivative and compare it to the continuous fractional derivative via a series of numerical experiments.
%
\end{abstract}

\end{frontmatter}


\section{Introduction}

Discrete exterior calculus (DEC) is a computational toolkit based on the principle of creating discrete analogues of objects and operators from smooth differential geometry.
Since the original thesis work of Hirani~\cite{Hirani2003}, DEC has been used in computer graphics applications, such as texture mapping~\cite{Dominitz2010}, direction field design~\cite{Vaxman2016}, fluid simulation~\cite{Elcott2007}, and meshing~\cite{mullen2011hot};  
geometry processing applications, such as symmetry detection~\cite{ben2010discrete}, spectral mesh processing~\cite{zhang2010spectral}, mesh parameterization~\cite{hormann2007mesh}, and matching~\cite{zaharescu2009surface}; 
and computational simulation, such as geometric mechanics~\cite{leok2004foundations}, Lie advection~\cite{mullen2011discrete}, stress fields in continuum mechanics~\cite{kanso2007geometric}, Hamiltonian PDEs~\cite{bridges2006numerical}, and incompressible fluids~\cite{Elcott2007}.

Absent from these applications is a consideration of how to apply the theory in the context of \textit{fractional partial differential equations} (FDEs).
A PDE is called ``fractional'' if it involves a non-local interaction of some kind, such as a time-dependent problem with ``memory'' or a spatial variable whose value at a point depends on values within some radius $\epsilon$ that cannot be taken arbitrarily small.
Such PDEs have terms involving \textit{fractional derivative operators} $({\partial}/{\partial x})^s$, where $s$ is a non-integer number; the definition of these operators will be discussed later.
Some recent examples of FDE-based modeling include:
\begin{itemize}
\item Roop~\cite{Roop2006} studied models of ``anomalous diffusion,'' in which a fractional advection-dispersion equation gave a more accurate description of fluid flow through porous media than a standard diffusion equation;
\item Ozgen et al.~\cite{Oktar2014} simulated cloth movement underwater using an FDE with a drag term, which gave visibly better results than standard techniques; 
\item Farquhar et al.~\cite{Farquhar2018} modeled electric currents through a heart with tissue affected by ischaemia using a mix of standard and fractional versions of the monodomain equation; they vary the fractional order across the computational domain according to the local level of blood supply.
\item Paquet and Viktor~\cite{PV2013} studied shape analysis by using a fractional deRham operator.  They used an eigenvalue decomposition to compute the fractional power of their operator, and found that the nonlocal nature of the fractional operator better allowed them to capture the shape of an object.
\end{itemize}

While many additional applications of FDEs can be found in the literature, an accurate, generalized, and efficient approach to their discretization remains elusive.
Our long term goal is to provide such an approach via the tools and framework provided by discrete exterior calculus.  As a starting point, in this paper, we introduce certain fundamental operators that would come up in a discrete approximation of a solution to a FDE.  Specifically, we focus on techniques for discretizing fractional derivative operators in accordance with DEC principles and conventions.

\begin{figure}[h!]
 \centering
    \begin{subfigure}[b]{0.14\textwidth}
        \includegraphics[width=\textwidth]{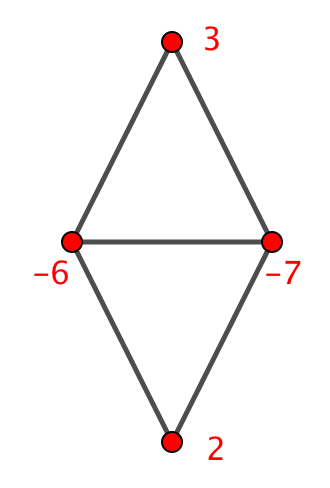}
        \caption{} \label{fig:0-form}
    \end{subfigure}
    ~ 
    \begin{subfigure}[b]{0.14\textwidth}
        \includegraphics[width=\textwidth]{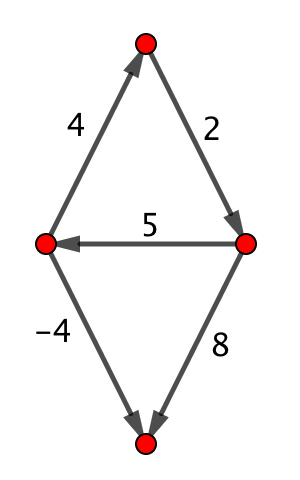}
        \caption{} \label{fig:1form}
    \end{subfigure}
    ~ 
    \begin{subfigure}[b]{0.14\textwidth}
        \includegraphics[width=\textwidth]{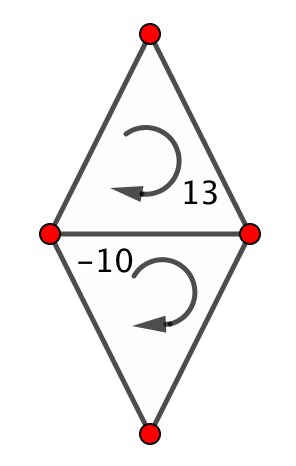}
        \caption{}\label{fig:2-form}
    \end{subfigure}
	\caption{A discrete 0-form, 1-form, and 2-form, (a)--(c), respectively.  A discrete fractional exterior derivative operator should incorporate ``nearby'' values when computing, more than just the incident values computed by the DEC derivative.}
\label{fig:discforms}	
\end{figure}

A brief review of key concepts from smooth and discrete exterior calculus helps motivate our approach.
In continuous exterior calculus, $p$-forms are used to generalize the description of scalar fields (0-forms), vector fields (1-forms), and tensor fields ($p$-forms, $p\geq 2$) defined over a smooth manifold (see e.g.~\cite{AMR2012} for a formal treatment).  
The continuous exterior derivative operator takes $p$-forms to ($p+1$)-forms, generalizing the notions of grad, curl, and div from vector calculus.
In DEC, discrete $p$-forms are defined over a simplicial complex, with discrete $p$-forms represented by values assigned to each $p$-simplex in the complex; see Figure~\ref{fig:discforms}.
The discrete exterior derivative operator $\mathbb{D}_p$ takes a discrete $p$-form to a discrete ($p+1$)-form by summing the values on the boundary of a $p+1$ simplex with signs $\pm1$ according to an orientation convention.

For example, the $\mathbb{D}_0$ operator that could be applied to the 0-form in Figure \ref{fig:0-form} is a (\# of edges)$\times$(\# of vertices)=$5\times 4$ matrix, where the entries are either $0$ or $\pm 1$ according to edge-vertex incidence:

\begin{align*}
\mathbb{D}_0 &= \begin{bmatrix}
      -1 & 1 & 0 & 0\\
  	1 & 0 & -1 & 0\\
  	0 & -1 & 1 & 0\\
  	0 & -1 & 0 & 1\\
  	0 & 0 & -1 & 1
\end{bmatrix}
.\end{align*}

We observe that a discrete fractional derivative operator in this context should still take discrete $p$-forms to discrete ($p+1$)-forms, but must somehow incorporate more than just vertex incidence information, according to the weight $s$ of the fractional derivative that is involved.  
Toward that goal, we present the following results:
\begin{itemize}
    \item Summarize relevant background from DEC and fractional calculus (Section~\ref{sec:bkgd}).
	\item Give a computable definition for the fractional derivative that adheres to the structure laid out in DEC (Section~\ref{sec:def-frac-deriv}) and investigate its properties (Section~\ref{sec:fd-props}).
	\item Test our definition on 1D and 2D examples (Section~\ref{sec:num-exp}) and discuss its effectiveness and possible future applications (Section~\ref{sec:disc}).
\end{itemize}

\section{Background and Notation}
\label{sec:bkgd}

In this section, we will give some background information to both DEC and fractional calculus.
\subsection{Discrete Exterior Calculus}
To use DEC, we build up a simplicial complex $K$ of dimension $N$ by using $p$-dimensional simplices $\sigma^p$, where $p \in \{0, 1, \ldots, N\}$.  These simplices are defined by $\sigma^p = [v_0, v_1 \ldots, v_p]$, where the order that we give $v_0, v_1, \ldots v_p$ specifies the orientation of the $p$-simplex.

On each $p$-simplex, we can attribute some value, and this represents a discrete $p$-form, also called a cochain.  To move from $p$-forms to $(p+1)$-forms, we use the discrete exterior derivative $\mathbb{D}_p$, which is a rectangular matrix that has a number of columns equal to the number of $p$-simplices in the complex and a number of rows equal to the number of $(p+1)$-simplices in the complex.  Specifically, $\mathbb{D}_p$ is the transpose of the boundary operator matrix.

The discrete exterior derivative is a local operator that does not make use of any metric information.  The matrix $\mathbb{D}_p$ only has nonzero entries of $\pm 1$ for simplices that are adjacent.  However, other DEC operators (such as the Hodge star) make use of metric information that is given by a local metric $d(v_0,v_1)$ where $v_0, v_1$ are part of an edge $[v_0, v_1] \in K$ (see \cite{Hirani2003} a discussion on this).  These values give information about the distances between vertices on the primal mesh.  As we will see later, this will become important when extending the discrete exterior derivative to a fractional setting.

\subsection{Fractional Calculus}

\begin{figure}[!ht]
\centering
	\includegraphics[width=0.9\columnwidth]{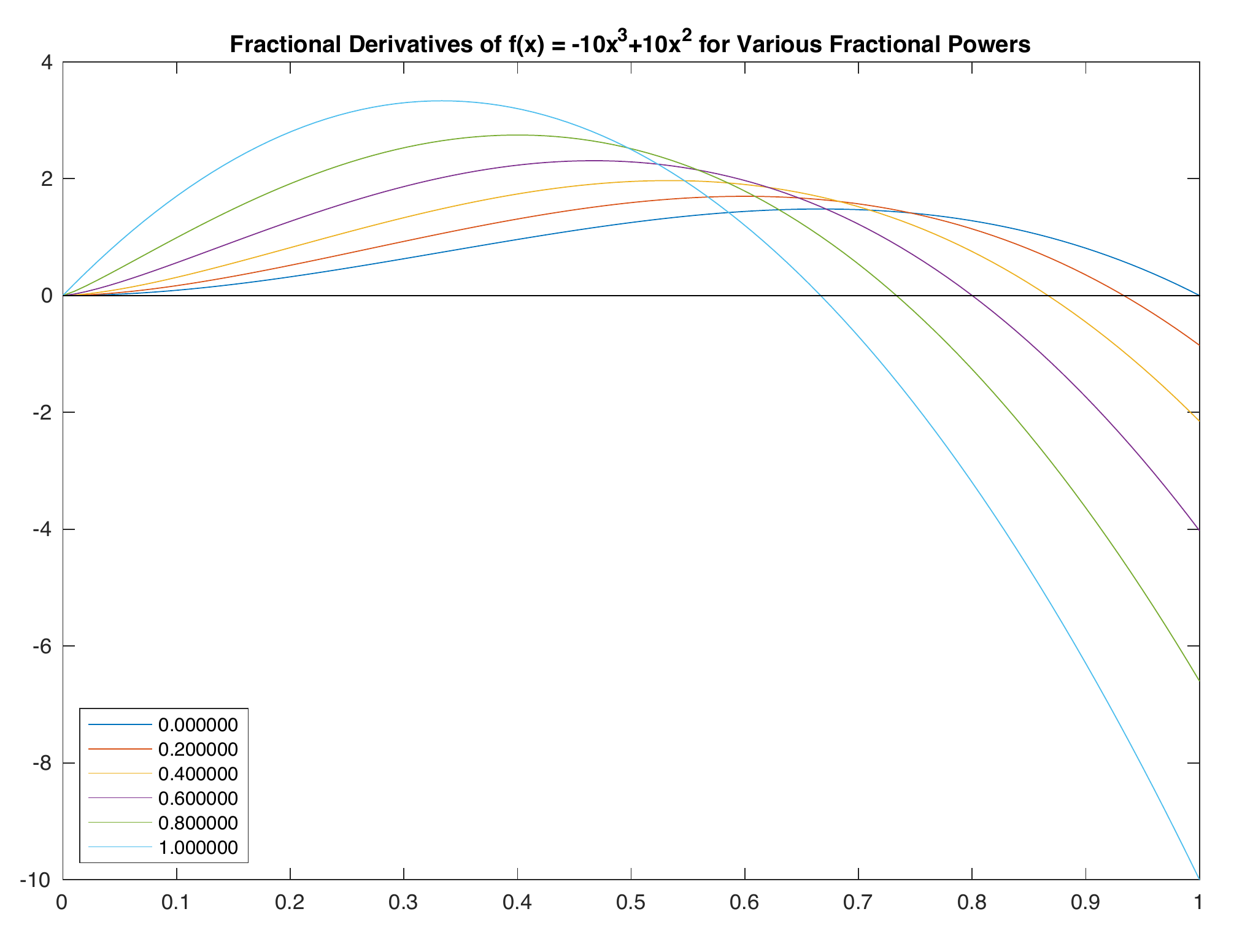}
        \caption{Fractional powers of order $s=0,~0.2,~0.4,~0.6,~0.8,~1$, for the left-sided Caputo derivative of $f(x) = -10x^3 + 10x^2$.} \label{fig:ExFractionalDerivatives}
\end{figure}

Fractional calculus has many different definitions of a fractional derivative (e.g. Riemann-Liouville, Caputo, Grunwald-Letnikov~\cite{Herrmann2010}).    
Central to all of them is the idea that the ``$s$th'' derivative of a function for $s \in [0,1]$ should give back a function in between the 0th (identity) and 1st derivative, varying continuously with $s$; see \autoref{fig:ExFractionalDerivatives}.
In our work, we will focus on the Caputo derivative.  Before we give the definition of the Caputo derivative itself, first we give some background information. 

\begin{definition}
	A function $f$ is said to be \emph{$n$-times continuously differentiable} if we can take $n$ derivatives of $f$ and $f^{(n)}$ is a continuous function.  In this case, we denote it as $f \in C^n[a,b]$.
\end{definition}

\begin{definition}
	The \emph{Gamma function} is given by
	
	\begin{equation}
		\Gamma(z) = \int_0^\infty e^{-t} t^{z-1} dt
	\end{equation}
	
	where $Re(z) > 0$.
\end{definition}
For natural numbers, the Gamma function acts like a factorial.  That is, $\Gamma(n+1) = n!$ whenever $n \in \mathbb{N}$.  For our purposes, the Gamma function will typically just be a constant that we can compute when needed.
Now we have the fractional Caputo derivative.

\begin{definition}
	The \emph{left-sided fractional Caputo derivative} of order $s \in \mathbb{R}^+$ of a function $f \in C^n[a,b]$ is 
	
	\begin{equation}
		D^{s}_{[a,x]} f(x) =
		\begin{cases}
			\displaystyle\frac{1}{\Gamma( n-s )} \int ^x _a \frac{f^{(n)}(\tau) d\tau}{(x-\tau)^{s +1 - n}} & s \notin \mathbb{N}\\[4mm]
			f^{(s)}(x) & s \in \mathbb{N}
		\end{cases} 
	\end{equation}
  where  $n = \ceil s$ (see \cite{Caputo1966}).  Note that often in the fractional calculus literature, the Caputo derivative is denoted as
 	
 	\[^C_aD_x^s f(x).\]
\end{definition}
 	
\noindent We choose to not use this notation because the majority of this work will focus only on the Caputo derivative.  The only time we deviate from this is in Section~\ref{sec:disc}, where we will briefly discuss the Riemann-Liouville derivative.  In that case, we will make sure it is clear that we are discussing a fractional derivative operator that is not the Caputo derivative.  

The subscript $[a,x]$ serves two purposes.  The first is that it defines the domain of integration for the fractional derivative, and the second is that it is used to denote which variable this is with respect to. 
We will also make use of a 2-sided definition of a fractional Caputo derivative.

\begin{definition}

   A \textit{two-sided fractional Caputo derivative} of order $s \in \mathbb{R}^+ \backslash \mathbb{Z}$ of a function $f \in C^n[a,b]$ is
    
    \begin{equation}
    D^s f(x) = \frac{1}{\Gamma(n - s)} \left(\int_a^x \frac{f^{(n)}(\tau) }{(x-\tau)^{s + 1 - n}}d\tau + \int_x^b \frac{f^{(n)}(\tau) }{(\tau - x)^{s + 1 - n}}d\tau \right).
    \end{equation}
\end{definition}

This can be written in a more compact form as 

\begin{equation}
D^s f(x) = \frac{1}{\Gamma(n - s)} \int_{[a,b] \backslash \{x\}} \frac{f^{(n)}(\tau) }{\vert x-\tau \vert ^{s + 1 - n}}d\tau.
\end{equation}

To understand the behavior of these operators, we have some examples in the continuous setting of the left-sided fractional Caputo derivative.
\begin{example}
	Due to taking a derivative before integrating, any constant function $f(x) = c$ on $[0,1]$ has the property that
	
	\begin{equation}
		D^{s}_{[a,x]} f(x) = 0
	\end{equation}
	
\noindent for any $s \in \mathbb{R}$.  
\end{example}

\begin{example}
	Let $f(x) = x^q$, $q > 0$ and $s \in (0,1)$.  Then we have that
	\begin{equation}
	D^{s}_{[0,x]} f(x) = \dfrac{\Gamma(q+1)x^{q-s}}{\Gamma(q+1-s)}.
		\label{eqn:CapDerPower}
	\end{equation}
	
\noindent This result mimics what we would expect in the case of $s = 1$, i.e.~$D^{1}_{[0,x]}x^q=qx^{q-1}$.  We see some results for the left-sided fractional Caputo derivative for $f(x) = -10x^3 + 10x^2$ in \autoref{fig:ExFractionalDerivatives} for a few different powers of $s$.
\end{example}

By changing the subscript notation to $[a_x, x]$ (or more generally $[a_i, x_i]$), we can extend the one-variable definition of the Caputo derivative to functions of multiple variables by writing

\[ {\frac{\partial ^{s}}{(\partial x)^{s}}}\bigg\rvert_{(a_x,x)} f(x,y) \]

\noindent with the understanding that this is just 
 
 \begin{equation}
 	\frac{1}{\Gamma(1-s)} \int_{a_x}^x \frac{f_x(\tau,y) d\tau}{(x-\tau)^{s}}
 \end{equation}
 
\noindent if $s \in (0,1)$.  
 
Regardless of how many variables we choose to work in, part of the difficulty with FDEs becomes obvious from the definitions: the operators that we are using are non-local.  Since these derivatives require the integration over some domain $[a,x]$, a small neighborhood of information at a point is not enough to give the value of the fractional derivative at that point.    This leads to some numerical difficulties.  



Primarily, a major issue is that numeric computations are more challenging.  For spatial FDEs, examining a large neighborhood manifests as more computational work to gather and process the surrounding information.  In general, for every fractional exponent, we may need to examine the values defined on the entire domain to accurately compute the derivative.  For fractional derivatives that involve time, the computation necessitates storing historical information for all previous values.  This represents a process with memory, but there is a tradeoff in that storing all history in practice would be computationally infeasible.

In our paper we will be focusing on spatial FDEs, but this manifests as its own set of computational problems.  Consider using a matrix to encode the derivative operator.  Due to the non-local nature of a fractional derivative, such matrices lose a lot of the sparsity and structure that they would normally have.  For example, in \cite{Roop2006}, Roop investigated the use of a finite element scheme to solve the fractional advection-dispersion equation.  In order to use quadratic basis elements on a $64 \times 64$ mesh, he had to store a full $16641 \times 16641$ stiffness matrix with no special structure.

Spatial FDEs are a primary factor in choosing to use a 2-sided fractional Caputo operator.  Using a left-sided operator in a temporal setting can be interpreted as meaning the process has a memory effect.  A left-sided operator in a spatial setting would ignore non-local effects that occur on the right side with no geometrical justification.  Thus the 2-sided definition we have given will be used in our discretization for defining a fractional discrete exterior derivative.

We will see a similar issue of having to store large, dense matrices later in this work.  Our main goal at this stage is the definition of a fractional discrete exterior derivative. Addressing how to effectively and accurately do this computations without losing the information the operator is trying to provide will be considered in future work.



\section{Defining a fractional discrete exterior derivative}
\label{sec:def-frac-deriv}

Combining fractional calculus and DEC comes with its own set of challenges.  Part of the design of DEC is that it is a set of tools that work in a local setting, with no need to understand global distances.

Our goal is to create a definition that imitates the fractional Caputo derivative.  While doing this, our definition for fractional discrete exterior derivative should also act in a similar fashion to the normal discrete exterior derivative map.  Thus we have some properties that we seek to preserve:

\begin{enumerate}
	\item That the operator $\mathbb{D}_{0}^s$ maps a discrete $0$-form $f$ with values at the vertices of a mesh to a discrete $1$-form $\mathbb{D}_{0}^s f$, taking values at the edges of the mesh.  Similarly, $\mathbb{D}_{1}^s$ should map a discrete $1$-form to a discrete $2$-form taking values at the faces.
	\item $\mathbb{D}_{0}^s$ maps a constant $0$-form to $0$ at each edge of the mesh.
	\item We want that as the fractional order goes to $1$, we recover $\mathbb{D}_k$.
\end{enumerate}

\noindent Note: the Caputo derivative does not yield the identity operator in the limit $s \rightarrow 0$, so it is not a property that we try to emulate. 

We mention this here before illustrating the design of our definition.  In choosing to use DEC to tackle fractional differential equations, we must make certain design choices that are going to lead to complicating portions of DEC that are normally straightforward.  We believe that the extra complications are worth it to make use of the framework that DEC gives us in solving differential equations.  In particular, the ability to computing distances between $p$-simplices is a non-trivial requirement for the computation.   We will discuss this more in section 4.

\subsection{The Fractional Discrete Exterior Derivative}
As we build our definition, we will start with the Caputo derivative definition and replace smooth operators with DEC operators wherever possible.  For convenience, we restrict ourselves to the setting where the two-sided fractional Caputo derivative is:

\[ D^s f(x) = \frac{1}{\Gamma(1 - s)} \int_{[a,b] \backslash \{x\}} \frac{f'(\tau) d\tau}{\vert x-\tau \vert ^s}. \]

\noindent whenever $s \in (0,1)$.

Before we begin discretizing, notice that the Caputo derivative operator can be decomposed into a  sequence of standard operations working inside out.  First we differentiate $f$, then we integrate a weighted version of the derivative.  Thus we will end up discretizing this in the same sequence of steps.

To discretize, we first apply a discrete exterior derivative to a $0$-form $\alpha$, yielding a discrete $1$-form $\mathbb{D}_0 \alpha$. The discrete $1$-form should take values at each edge, so for the remaining steps we will focus on what the fractional operator computes at a specific edge $x$. 

Next, we discretize the integral with a kernel of $\dfrac{1}{\vert x-t \vert^s}$.  This is similar to discretizing a convolution of $\mathbb{D}_0 \alpha$ against the convolution kernel of $\dfrac{1}{\vert x-t \vert^s}$.   This translates to computing a weighted sum 
\[ \sum_{i=0}^{\vert E \vert -1} \dfrac{1}{\vert \vert x - t_i \vert \vert ^s} \mathbb{D}_0 \alpha\]
where $t_i$ runs over each of the edges in the mesh, $x$ is the edge of interest, and $\vert \vert x - t_i \vert \vert$ represents the distance between the barycenters of the two edges.

However, this is undefined when $t_i = x$.  To fix this, we will denote a weighting vector 
\begin{equation}
	(w_{x})_i = \dfrac{1}{\vert \vert x - t_i \vert \vert ^s }
	\label{eqn:weights}
\end{equation}
whenever $ x \neq t_i$
 and 
\begin{equation}
 	(w_{x})_i = C_s \text{max}_{j \neq k} \left( \dfrac{1}{ \vert \vert t_k - t_j \vert \vert ^s}\right)
\end{equation}
when $x = t_i$.  Note that the constant $C_s$ is used to give more weight to the current edge than the other edges in the mesh.  Choosing how to properly define this entry in the vector $w_x$ is an interesting problem by itself, and our definition is sensitive to the choice of this value.  We found in our numerical experiments that that $C_s=\frac{2s}{1-s}$ is a reasonable choice. Finally, we weight by $\frac{1}{\Gamma(1-s)}$.  This yields the value of $(\mathbb{D}_0^s \alpha)(x)$. 
\[ 
\left(\mathbb{D}_0^s \alpha \right)(x) = \dfrac{1}{\Gamma(1-s)}  w_{x} \cdot  \mathbb{D}_0 \alpha .
\]

We can generalize this definition as follows:

\begin{definition}[Fractional Discrete Exterior Derivative]
	\label{def:FDED1}
	Let $\alpha$ be a discrete $p$-form on a connected simplicial complex $M$.  Then the \emph{$(p,s)$-fractional discrete derivative of $\alpha$} at a $(p+1)$ simplex $\sigma^{p+1}$ is 
	\begin{equation}
	 	\left(\mathbb{D}_{p}^s \alpha\right)(\sigma^{p+1}) = \dfrac{1}{\Gamma(1-s)} w_{\sigma^{p+1}} \cdot \mathbb{D}_p \alpha
	\end{equation}
where $\sigma^{p+1}$ is a specific $p+1$ simplex, $w_{\sigma^{p+1}}$ is the weighting vector from above, replacing edges $x_k$ and $t_i$ with $(p+1)$-simplices as necessary, and $\mathbb{D}_p \alpha$ is the typical $p^{th}$ order discrete exterior derivative of $\alpha$.  
\end{definition}

\section{Properties of the fractional discrete exterior derivative}
\label{sec:fd-props}


Before we get to testing the definition on some examples, we proceed with a brief discussion on the properties of the above \autoref{def:FDED1}.  In this discussion, $\vert P \vert$ will represent the number of $p$-simplices in the mesh.

\begin{enumerate}
	\item This definition does enact the property that $\mathbb{D}_p^s$ correctly maps from $p$-forms to $(p+1)$-forms. To see this explicitly, see the equivalent matrix form of the definition that we give below.
	\item  We know from DEC that the dot product
\[  w_{\sigma^{p+1}} \cdot \mathbb{D}_p \alpha \]
	will result in $\mathbb{D}_p^s \alpha = 0$ whenever $\alpha$ takes a constant value for all $p$-simplices in the mesh. 

	\item Currently, this definition does \emph{not} yield a recovery of $\mathbb{D}_p$ when we take $s \rightarrow 1$.  This can be fixed, and we will show what this may look like below.
\end{enumerate}

Note that we can rewrite the definition of the fractional discrete exterior derivative operator using matrices.  To do this, let $W$ be the matrix made out of using the vectors $w_{x_k}$ as the columns. This creates a $\vert P \vert \times \vert P \vert$ matrix that we can use to say 
	\begin{equation}
		\mathbb{D}_p^s \alpha = \frac{1}{\Gamma(1-s)} W\mathbb{D}_n\alpha. 
	\end{equation}

\noindent The choice to start this way was made to make it easier to see the action at a specific simplex.  However, this form makes it clear that $\mathbb{D}_p^s : \Omega_d^p \rightarrow \Omega_d^{p+1}$.





The natural way that we can force the property of recovering $\mathbb{D}_p$ in the limit of $s \rightarrow 1$ is to piecewise define $w_{x_k}$ based off the fractional order.  If $s \in (0,1)$ (or more generally, $s \in \mathbb{R}^+ \backslash \mathbb{N}$), then define $w_{x_k}$ as in \autoref{eqn:weights}, and in the case that $s = 1$ (or $s \in \mathbb{N}$), define
 \begin{equation} 
 	(w_{x_k})_i = \delta_{ik},
 \end{equation}
 where $\delta_{ik}$ is the Kronecker delta function, taking a value of $0$ if $i \neq k$ and $1$ when $i = k$.  One reason for doing this piecewise definition is that $\Gamma(0) = 0$.  Therefore both the fractional Caputo derivative and the fractional discrete exterior derivative cannot be defined here by the usual definition or the new definition, respectively.  Furthermore, in order for $\mathbb{D}_p^s$ to recover $\mathbb{D}_p$ in the limit $ s \rightarrow 1$, we would need that $W \rightarrow \Gamma \left(1-s\right) I$.  In section 5, we give some evidence that as $s \rightarrow 1$, our approximation gets better.  This indicates that the $W$ matrix is indeed becoming diagonally dominant so that the fractional contribution gets less as we get close to $s = 1$.  Hence we believe that we do achieve property $3$.

There are other properties that we chose not to enforce that may have been natural given the setting of DEC.  For example, the composition property $\mathbb{D}_{p+1} \circ \mathbb{D}_p = 0$.  However, our definition does not satisfy an extension to $\mathbb{D}_{p+1}^s \circ \mathbb{D}_p^s = 0$.  To see this, let $W_k$ be the corresponding weighting matrix to $\mathbb{D}_k^s$ and notice that:
\begin{align}
	\mathbb{D}_{p+1}^s \mathbb{D}_p^s \alpha &= \frac{1}{\Gamma(1-s)} \left(W_{p+1} \mathbb{D}_{p+1}\right)\mathbb{D}_p^s \alpha \\
	&= \frac{1}{\Gamma(1-s)} W_{p+1} \mathbb{D}_{p+1} W_p \mathbb{D}_p \alpha. 
\end{align}

\noindent In general, $W_p$ and $\mathbb{D}_p$ cannot commute, and the distances that are involved in $W_k$ will not force $0$ from any specific entries.

Part of the expectation for $\mathbb{D}_{p+1} \circ \mathbb{D}_p = 0$ comes from the smooth setting, where we have $d^{p+1} \circ d^p = 0$. In other words, all exact forms in the smooth setting are closed.  While a similar property would be good to have, we do not have a clear reason for expecting it to be true after the weighting process occurs.  For additional formal discussion on the smooth theory of fractional exterior calculus, see chapter 12 of \cite{Tarasov}.

Finally, we must address the use of a norm $\vert \vert \cdot \vert \vert$ in the definition.  The DEC framework only assumes the existence of a local metric~\cite{Hirani2003}, which can be used to measure distance only between adjacent simplices.  We next discuss how to extend the local metric to be able to compute the distance between two arbitrary $p$-simplices.

To do this, we first need a well-defined way to find the distance between two arbitrary vertices on our mesh.  We make use of the fact that between any adjacent vertices $u$ and $v$, we know the distance $d(u,v)$ where $d$ is the local metric on the mesh.  Then we can apply an all pairs-shortest path algorithm to find the minimum distance between each pair of vertices on the path.  This gives a distance between any pair of vertices, not just adjacent ones, which we denote $d_m(u,v)$, where the subscript $m$ is used to differentiate the minimum distance between any pair of vertices on the mesh from the local metric $d$.

We extend distance between $0$-simplices to build a distance between two $p$-simplices.  Let $\sigma^p$, $\eta^p$ be two simplices, $\{u_k\}$ and $\{v_k\}$ be the set of vertices attached to $\sigma^p$ and $\eta^p$ respectively, with $k = 0, 1, \ldots, p$.  Then we define the distance
\[\vert \vert \sigma^p - \eta^p \vert \vert\]
computing
\[ \min_{i,j} d_m(v_i, y_j) + l(\sigma^p) + l(\eta^p) \]
for all pairs of $i, j$ with $i, j = 0, 1, \ldots, p$, where $l(\sigma^p)$ and $l(\eta^p)$ is the length between the barycenter and the boundary of each of the respective simplices.

As an aside, if we know that the embedding of the simplicial complex $K$ uses a metric that is induced by Euclidean metric of $\mathbb{R}^N$, then it is easier to implement the weighting matrix by using the Euclidean distance between circumcenters or barycenters of the two simplices.  That is, if we know that there exists a global metric that correctly embeds the complex, we can use the global metric.

\section{Numerical experiments}
\label{sec:num-exp}
With our definition in hand, we test it in some controlled scenarios.  First though, we give a brief outline for how one would go about implementing the fractional discrete exterior derivative as defined above, assuming we already have a discrete $p$-form  $\alpha$.
\begin{enumerate}
	\item Create $\mathbb{D}_p$ matrix.
	\item Create matrix $W$.
	\item Apply the product of $W$ and $\mathbb{D}_p$ to $\alpha$.
	\item Multiply by the scaling factor $\displaystyle\frac{1}{\Gamma(1-s)}$.
\end{enumerate}
In step $2$ we will have to do a computation that finds the distance between any two given $(p+1)$-simplices before $W$ can be created.

Now we do an example.  Recall from \autoref{eqn:CapDerPower} the form of the Caputo derivative acting on power functions
\[
D^{s}_{[0,x]} x^q = \dfrac{\Gamma(q+1)x^{q-s}}{\Gamma(q+1-s)}.
\]
We used the left-sided fractional Caputo derivative to find this result, and assumed that $s \in (0,1)$ with a domain of interest of $[0,1]$.  However, we based our fractional discrete exterior derivative on the 2-sided fractional Caputo derivative.  In the following examples, we follow the convention that $D^s_{[0,x]}$ represents a left-sided Caputo derivative, while the notation $D^s$ represents a two-sided Caputo derivative.


\begin{example}
	 As a specific example, consider $x^3$ on $[0,1]$ and $s = .5$.  Then we have that 
	\begin{equation}	 
	 	D^{.5} x^3 = \dfrac{\Gamma(4)x^{2.5}}{\Gamma(3.5)} - \frac{3}{\Gamma(3.5)} (1-x)^{\frac{1}{2}} (.75 + x + 2x^2)
	 \end{equation}
	 as the analytic solution.
	 
\end{example}

\begin{figure}[h!]
\centering
	\includegraphics[width=0.9\columnwidth]{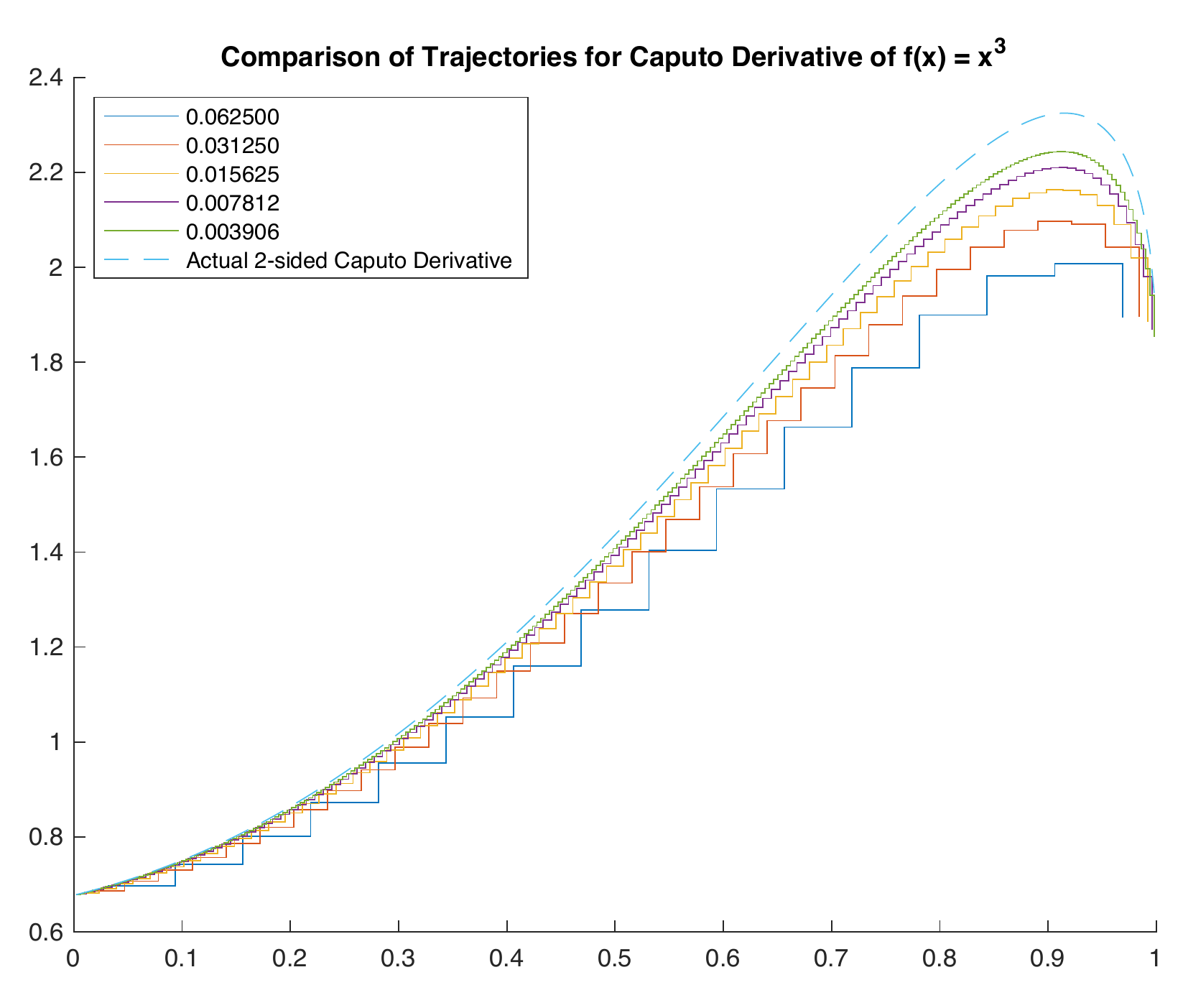}
        \caption{The comparison of $\mathbb{D}_0^{.5} x^3$ with $D^{.5} x^3$, for uniform meshes of various element widths.  The element widths tested are given in the legend.} \label{fig:Uncorrected1d}
\end{figure}

As we see in this example, according to \autoref{fig:Uncorrected1d} using our definition yields a similar values. The plots of the weighted discrete $1$-forms that result from the fractional discrete exterior derivative applied to $x^3$ on the mesh are obtained using the stairs plot in Matlab.  This plotting method results in a staircase function where a given step starts at a barycenter of an edge and lasts until the next step.  For each of these steps, we use the computed value in the vector $\mathbb{D}^s_0 x^3$, and plot it as a constant value for the rest of the edge.


\begin{example}
Let $s = 0.5$, $[a,b] = [0,1]$, and consider the function

	\[f(x) = -10x^3 + 10x^2, \]

with a fractional Caputo derivative of


\begin{equation}
\begin{split}
   D^{.5} f(x) = & -  \dfrac{10\Gamma(4)}{\Gamma\left(\dfrac{7}{2}\right)}x^{\frac{5}{2}} + \dfrac{10\Gamma(3)}{\Gamma\left(\dfrac{5}{2}\right)}x^{\frac{3}{2}} -  \\ & \dfrac{30}{\Gamma \left(\dfrac{7}{2}\right)}(1-x)^{\frac{1}{2}} \left(\frac{3}{4} + x + 2x^2\right) + \dfrac{20}{\Gamma \left(\dfrac{5}{2}\right)}(1-x)^{\frac{1}{2}}(x + \dfrac{1}{2}).
\end{split}
\end{equation}

\end{example}

Since the function we are working with has an easily calculated closed form for its fractional Caputo derivative and a discrete $1$-form can be thought of as piecewise constant on edges, we are able to symbolically integrate to get our error values.  The results of the error analysis are given in the table below.


\begin{center}
\begin{tabular}{ |c|c| }
	\hline
	Edges & $L_2$ error\\
	\hline
	2 & 1.5619\\
	\hline
	4 & 0.9933\\
	\hline
	8 & .6778\\
	\hline
	16 & .4759\\
	\hline
	32 & .3363\\
	\hline
	64 & .2378\\
	\hline
	128 & .1681\\
	\hline
	256 & .1188\\
	\hline
	512 & .0839\\
	\hline
	1024 & .0593\\
	\hline
\end{tabular}
\end{center}

\begin{figure}[h!]
\centering
	\includegraphics[width=1\columnwidth]{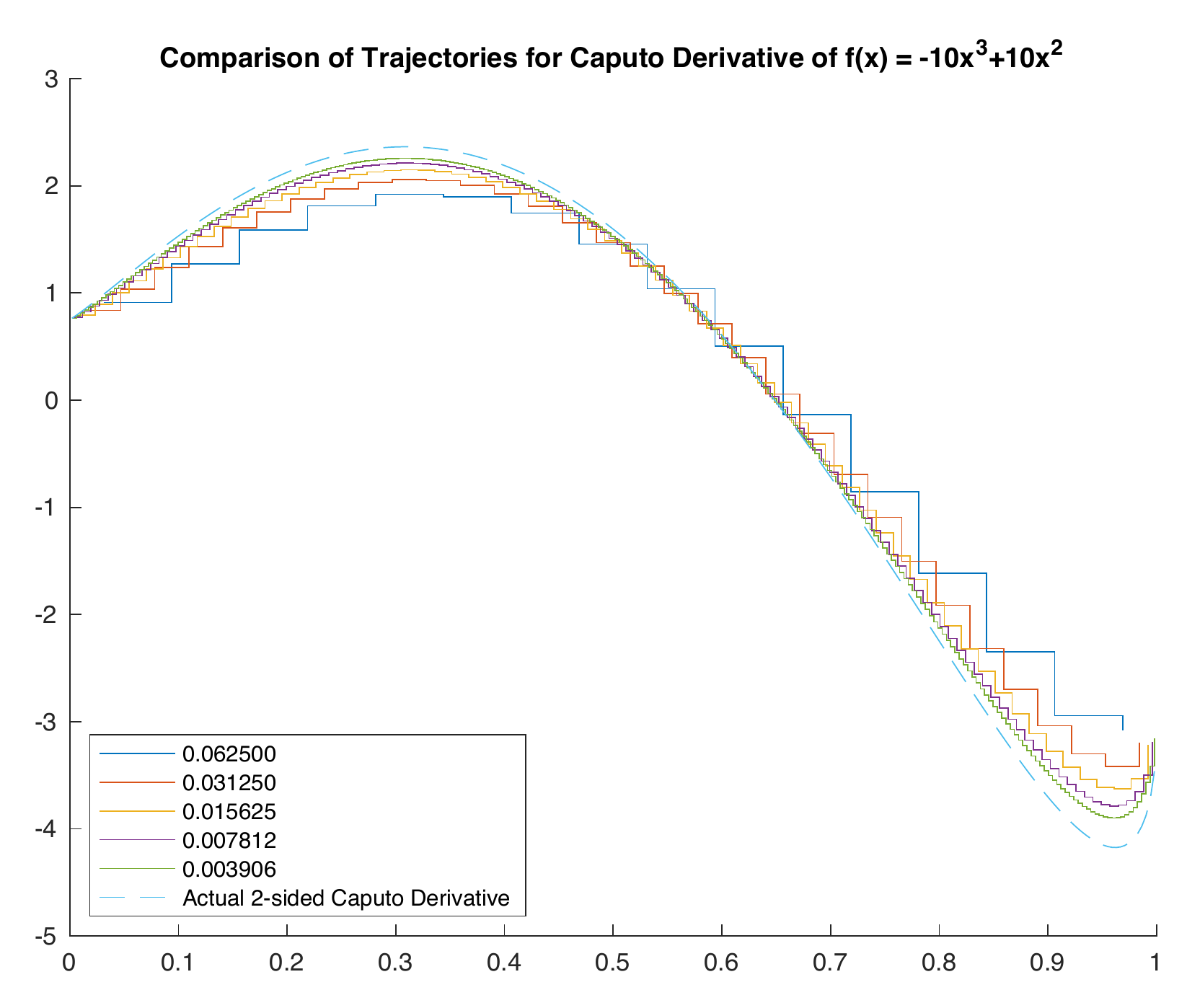}
        \caption{The comparison of $\mathbb{D}_0^{.5}\left( -10x^3+10x^2\right)$ with $D^{.5} \left(-10x^3+10x^2\right)$.  We use varying step sizes denoted in the legend.} \label{fig:1dEx2}
\end{figure}

The error results here give the indication that our error is going to $0$ as the step size decreases. In \autoref{fig:1dEx2}, we compare the closed-form version against another piecewise linear estimate that we get as the result of applying the fractional discrete exterior derivative.  In the case of different step sizes, we see that we do get closer to the true trajectory of the fractional derivative.

Since this example has a simple analytic solution to the fractional Caputo derivative we can test the effect of the value of $s \in (0,1)$. 

\begin{figure}[h!]
\centering
	\includegraphics[width=0.9\columnwidth]{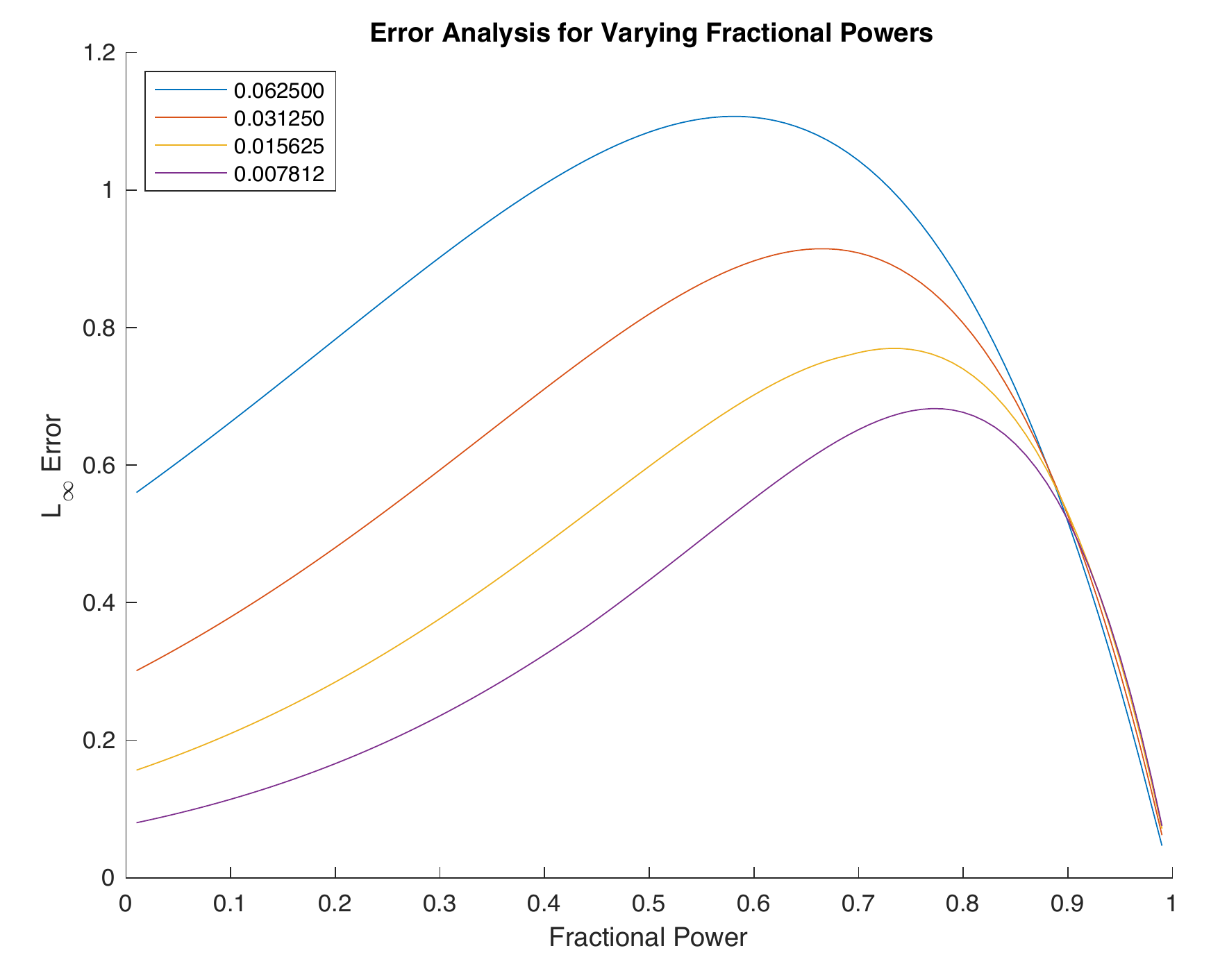}
        \caption{We compare the fractional exponent $s$ against the $L_\infty$ error for the function $f(x) = -10x^3 + 10x^2$.  We use a step size varying step sizes given in the legend.} \label{fig:VaryingFracError}
\end{figure}

In FDEs, it is commonplace to test the error against the fractional powers.  Some applications have a clear reasoning, whether experimental or physical, for desiring a fractional power to be some specific value.  Other applications treat the fractional power as free parameter.  Our results in \autoref{fig:VaryingFracError} show that while we can change the fractional power as we desire in $(0,1)$, our error can fluctuate within this range.  Part of this error is expected because of the $\frac{1}{\Gamma(1-s)}$  constant based on $s$.

Finally, we want to check how our definition compares to the left-sided fractional Caputo derivative of $e^x$.  To that end, we give the following definition and then a brief discussion of the modification of $W$ to account for a left-sided fractional derivative instead of a 2-sided fractional derivative.

\begin{definition}
The Mittag-Leffler two parameter function is
\begin{equation}
	E_{a,b}(z) = \sum_{k=0}^\infty \dfrac{z^k}{\Gamma(ak + b)}
\end{equation}
where $a, b \in \mathbb{R}^+$ and $z \in \mathbb{C}$.
\end{definition}
In general, we can think of this function as being a generalization of the exponential function.  In fact, for the values of $a = 1, b = 1$, we recover the exponential function exactly.  The Mittag-Leffler function can be used to give a nice looking form for the result of the fractional Caputo derivative of $e^x$. Explicitly, we get

\begin{equation}
	D_{[0,x]}^s e^{x} = x^{1-s} E_{1, 2 - s}(x)
\end{equation}

\noindent whenever $s \in (0,1)$ \cite{Mariya}.

\begin{figure}[h!]
\centering
	\includegraphics[width=0.9\columnwidth]{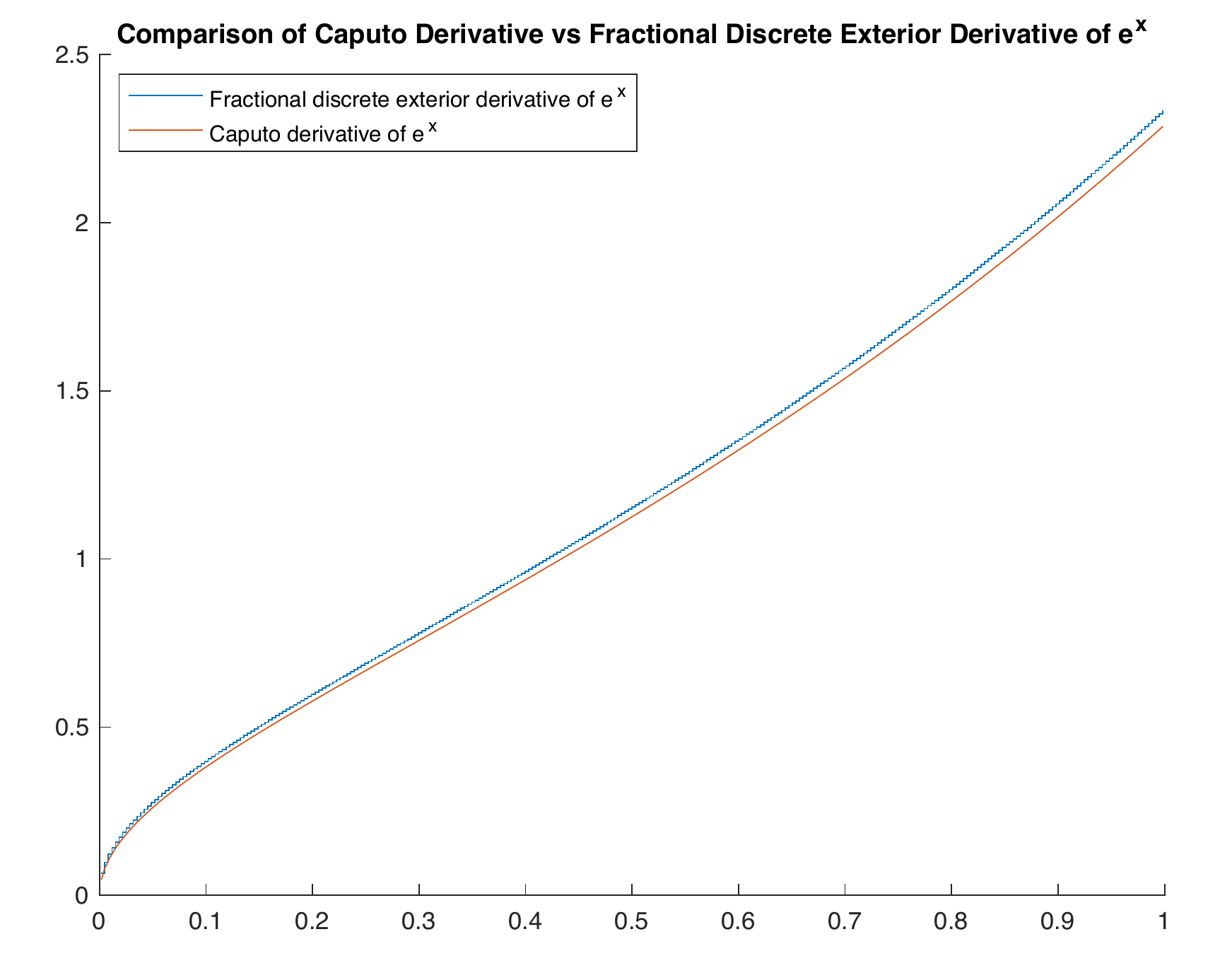}
        \caption{Comparison of the Caputo derivative of $e^x$ versus $\mathbb{D}_0^{.5} e^x$.} \label{fig:dex}
\end{figure}

Our discretization of the fractional discrete exterior derivative is based off the two-sided fractional Caputo derivative.  However, we want to be able to compare against $D_{[0,x]}^s e^{x}$.  In the left-sided fractional Caputo derivative, we have an integration domain of $[0,x]$.  To account for this, we need to modify the weighting matrix $W$.  Specifically, each entry of $W$ represents the weight assigned to one simplex in relation to another simplex, say $\sigma_i$ and $\sigma_j$.  Then, in the $1$-d case, to determine if the entry should be zero to account for the left-sided version of the derivative we compare the x-coordinates of the barycenters of $\sigma_i$ and $\sigma_j$.  If $x_i < x_j$, we leave $W_{i,j}$ as it is and otherwise we set $W_{i,j} = 0$.  The results of this comparison are given in \autoref{fig:dex}.

\begin{figure}[h!]
\centering
	\includegraphics[width=1\columnwidth]{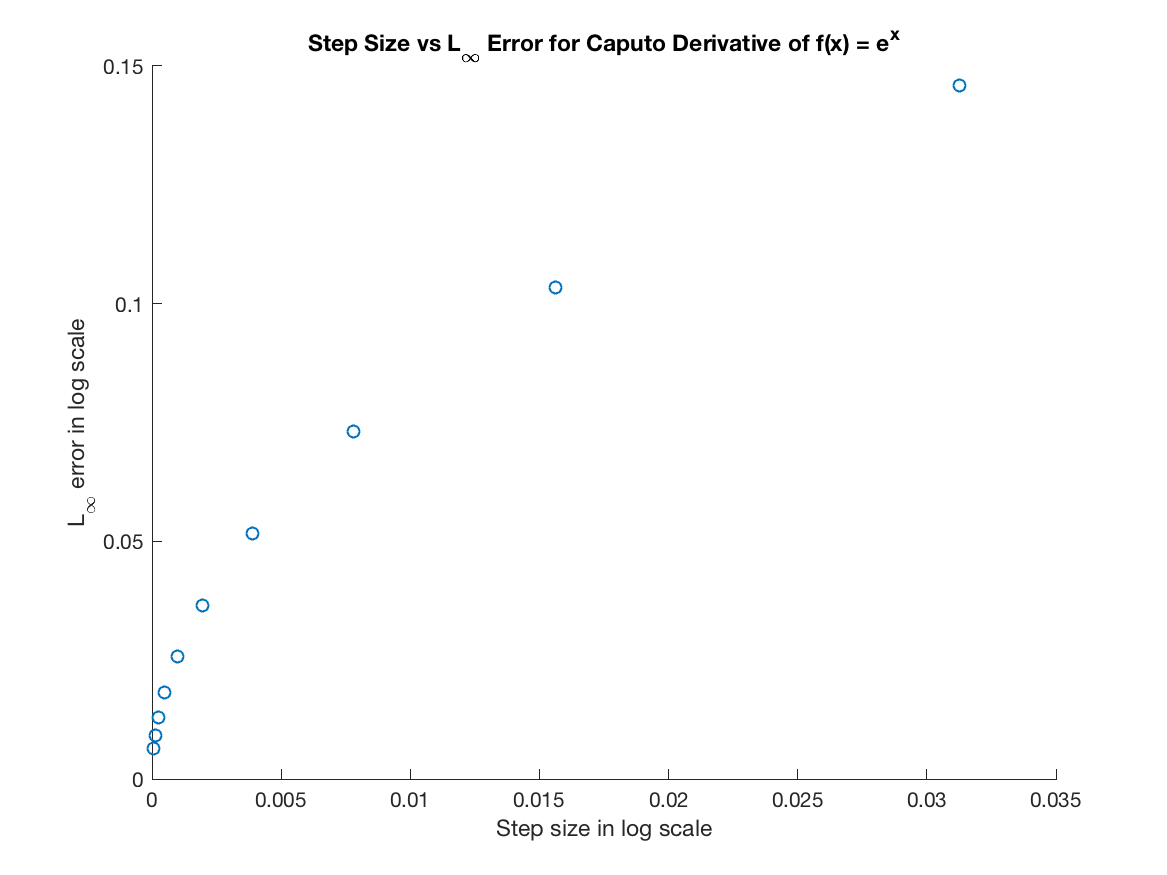}
        \caption{Error analysis for $\mathbb{D}_0^s e^x$.  We use the $L_\infty$ norm to find max errors at varying step sizes in comparison to $x^{.5} E_{1, 1.5}(x)$, the closed form of the left-sided Caputo derivative of $e^x$. } \label{fig:dexError}
\end{figure}

In \autoref{fig:dex}, we see the difference in trajectories for the Caputo derivative of $e^x$.  The code used to generate the Mittag-Leffler function is from \cite{Garrappa2015}.  Our definition undershoots the value of the Caputo derivative, but does seem to exhibit the correct behavior.  To do an accurate error analysis in this case, we choose the $L_\infty$ norm.  We see the results of this analysis in \autoref{fig:dexError}.

The following examples are efforts in testing the fractional discrete exterior derivative in $2$-d.  We start with a basic function

\begin{example}

Let $f(x,y) = -x^2 + y^2$.  Then the 2-sided fractional gradient field on the region $M = [0,1] \times [0,1]$ of order $s \in (0,1)$ is
\begin{equation}
	\left(\dfrac{-2x^\frac{3}{2}}{\Gamma \left(\frac{5}{2}\right)} - \dfrac{2(1-x)^\frac{1}{2}(x + \frac{1}{2})}{\Gamma \left(\frac{5}{2}\right)}, \dfrac{2y^{\frac{3}{2}}}{\Gamma \left(\frac{5}{2}\right)} + \dfrac{2(1-y)^\frac{1}{2}(y + \frac{1}{2})}{\Gamma \left(\frac{5}{2}\right)} \right).
\end{equation}

Using this, we consider a triangulated plot of $M$ with the gradient vector field drawn at the barycenters of the triangles.  This field is plotted in \autoref{fig:GradientFieldSaddle}.

\begin{figure}[h!]
\centering
	\includegraphics[width=0.9\columnwidth]{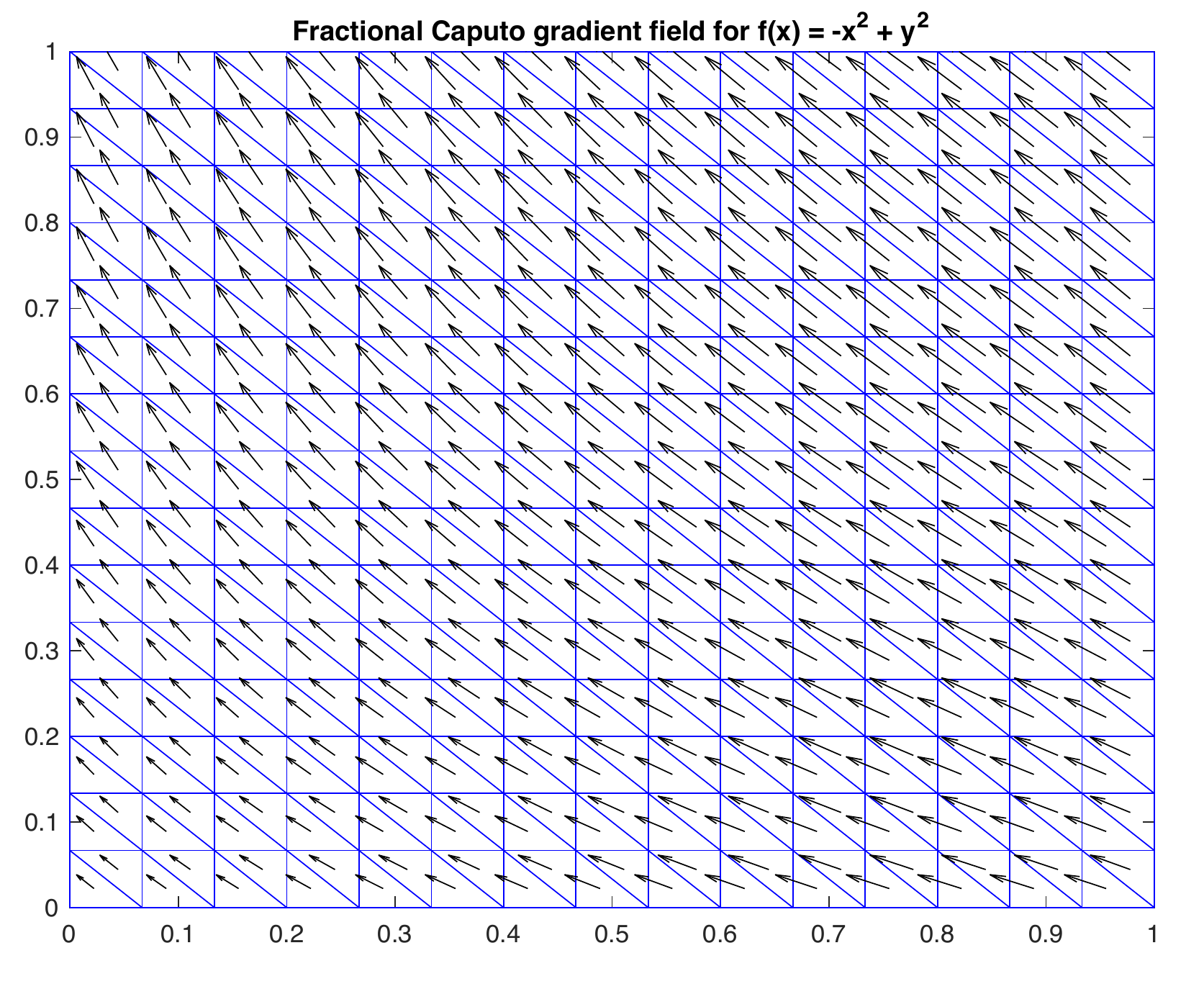}
        \caption{The 2-sided fractional Caputo gradient field for the function $-x^2 + y^2$. } \label{fig:GradientFieldSaddle}
\end{figure}

To test the fractional discrete exterior derivative on this example and visualize a vector field requires some extra work.  We first create the $0$-form vector $\alpha$ that represents the value of $f$ at the vertices of the mesh.  After applying the fractional discrete exterior derivative, we get a discrete $1$-form $\mathbb{D}^{.5} \alpha$, i.e. a cochain storing a value for each edge of the mesh with edge orientation implied by the global vertex ordering.  Visualizing and measuring the error of this data type requires some care. 

We first construct the Whitney map for 1-cochains, based on~\cite[Definition 3.3.4]{Hirani2003}. The construction is as follows: compute the barycentric coordinates $\lambda_1$, $\lambda_2$, $\lambda_3$ locally on each triangle, use these to compute the three Whitney 1-forms ($\lambda_i\nabla\lambda_j-\lambda_j\nabla\lambda_i$) associated to each edge, then weight each Whitney 1-form by the cochain value on the corresponding edge.  The result is a vector field, defined piecewise over the entire mesh, that recovers the cochain values on each edge (when projected onto an edge).  We evaluate this global vector field at the barycenter of each triangle for both qualitative and quantitative error analysis.

To give a sense of what the fractional Caputo gradient field on $M$ looks like, see \autoref{fig:GradientFieldSaddle}.  In \autoref{fig:RelErrorSaddle}, we look at relative error.  Notice that originally, the function has a saddle type critical point (i.e. the gradient of the function is zero) at the origin.   This critical point, along with some boundary effects from our definition, lead to a high relative error near the origin.  

\begin{figure}[h!]
\centering
	\includegraphics[width=0.9\columnwidth]{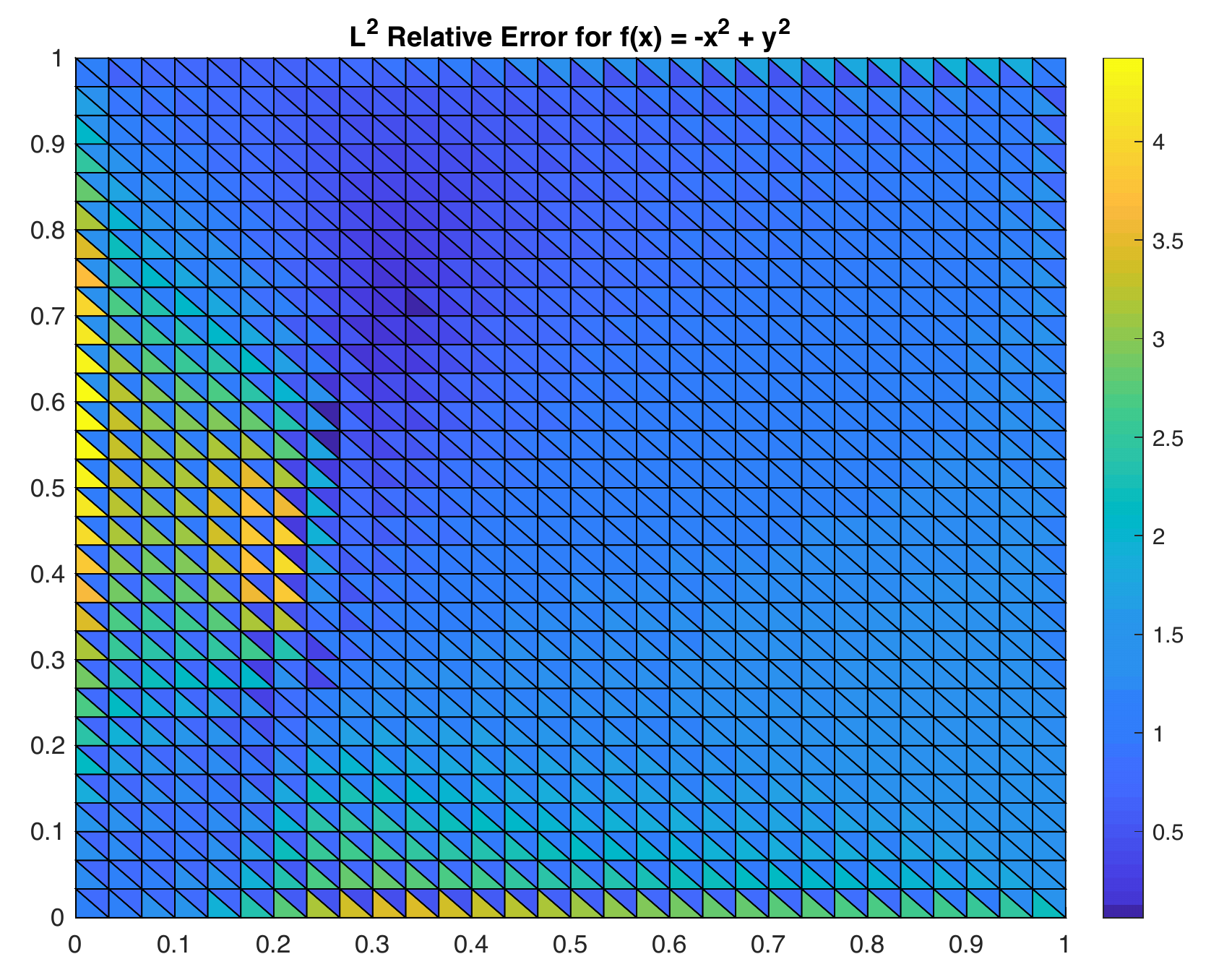}
        \caption{The triangles are colored to represent the relative $L^2$ error for $f(x,y) = -x^2+y^2$.  We get error values ranging between .0627 and 4.4227, with an average of 1.2126.} \label{fig:RelErrorSaddle}
\end{figure}
\end{example}

\begin{example}
We do one last example.  Let

\begin{equation}
	f(x,y) = (x-.1)^2 + (y-.1)^2
\end{equation}

be defined on $M = [0,1] \times [0,1]$.  In this case, we know that $f(x,y)$ has a minimum at $(.1,.1)$.  The 2-sided fractional Caputo gradient field of $f$ is 

	


\begin{equation}
\begin{split}
D^\frac{1}{2} f  = & \Bigg(\frac{2x^\frac{3}{2} - (1-x)^\frac{1}{2}(\frac{1}{2} - 2x)}{\Gamma \left(\frac{5}{2}\right)} - \frac{x^\frac{1}{2}-(1-x)^\frac{1}{2}}{5\Gamma \left(\frac{3}{2}\right)}, \\ 
& \frac{2y^\frac{3}{2} - (1-y)^\frac{1}{2}(\frac{1}{2} - 2y)}{\Gamma \left(\frac{5}{2}\right)} - \frac{y^\frac{1}{2}-(1-y)^\frac{1}{2}}{5\Gamma \left(\frac{3}{2}\right)} \Bigg)
\end{split}
\end{equation}

This fractional gradient field is plotted in \autoref{fig:GradientFieldMinimum}.  We also give a plot for relative error for $f$ in \autoref{fig:RelErrorMinimum}.  In this case, we again have a critical point that affects our results.  For this function, the critical point exists in the domain $M$.  One point of interest is that the fractional Caputo derivative does not preserve the location of critical points, which is what causes a large amount of the error in \autoref{fig:RelErrorMinimum}.  In either setting, understanding the behavior of the fractional discrete exterior derivative on functions of multiple variables with critical points is a topic that we intend to investigate further in the future.  

\begin{figure}[h!]
\centering
	\includegraphics[width=0.9\columnwidth]{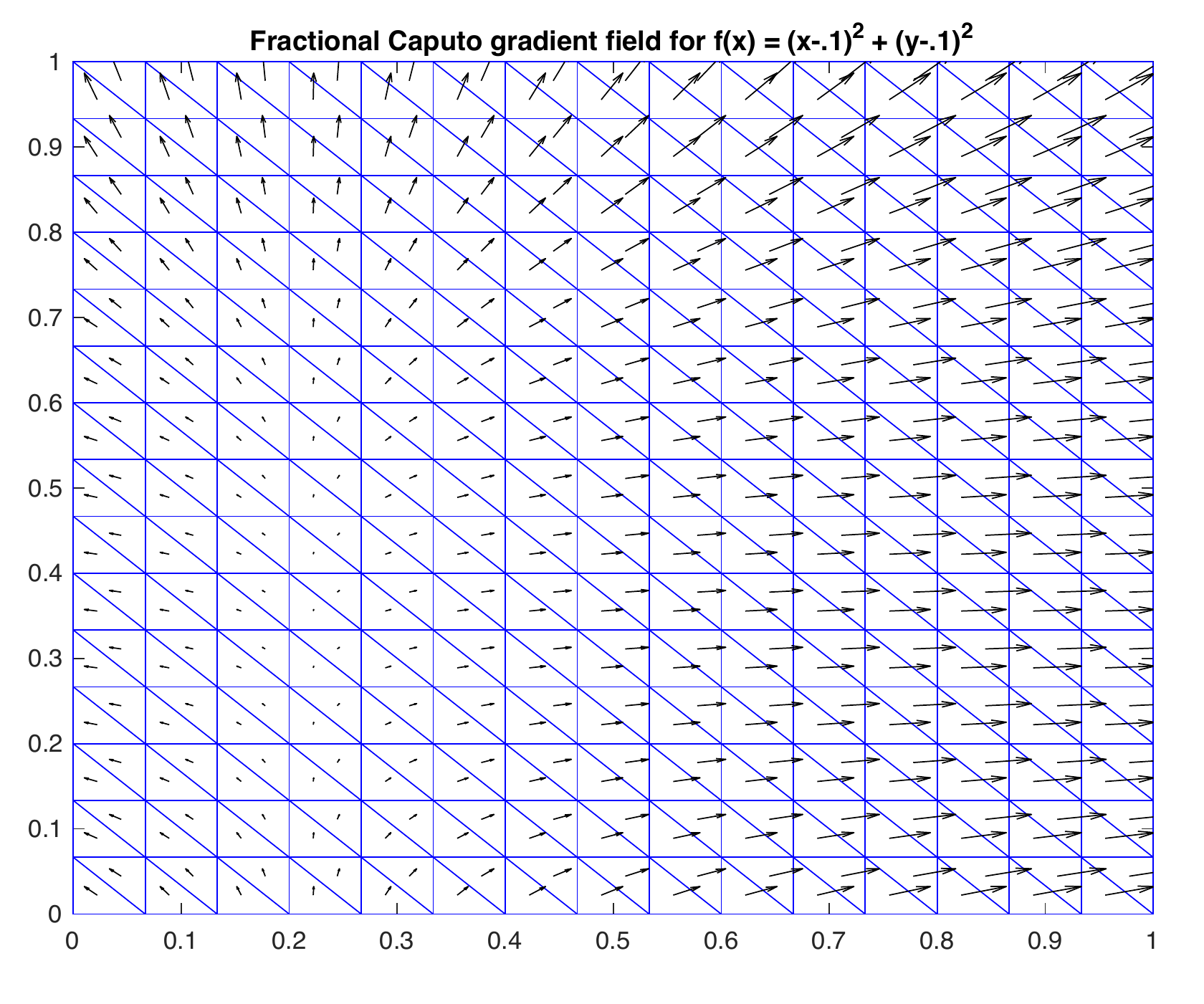}
        \caption{The 2-sided fractional Caputo gradient field for the function $f(x,y) = (x-.1)^2 + (y-.1)^2$. } \label{fig:GradientFieldMinimum}
\end{figure}

\begin{figure}[h!]
\centering
	\includegraphics[width=0.9\columnwidth]{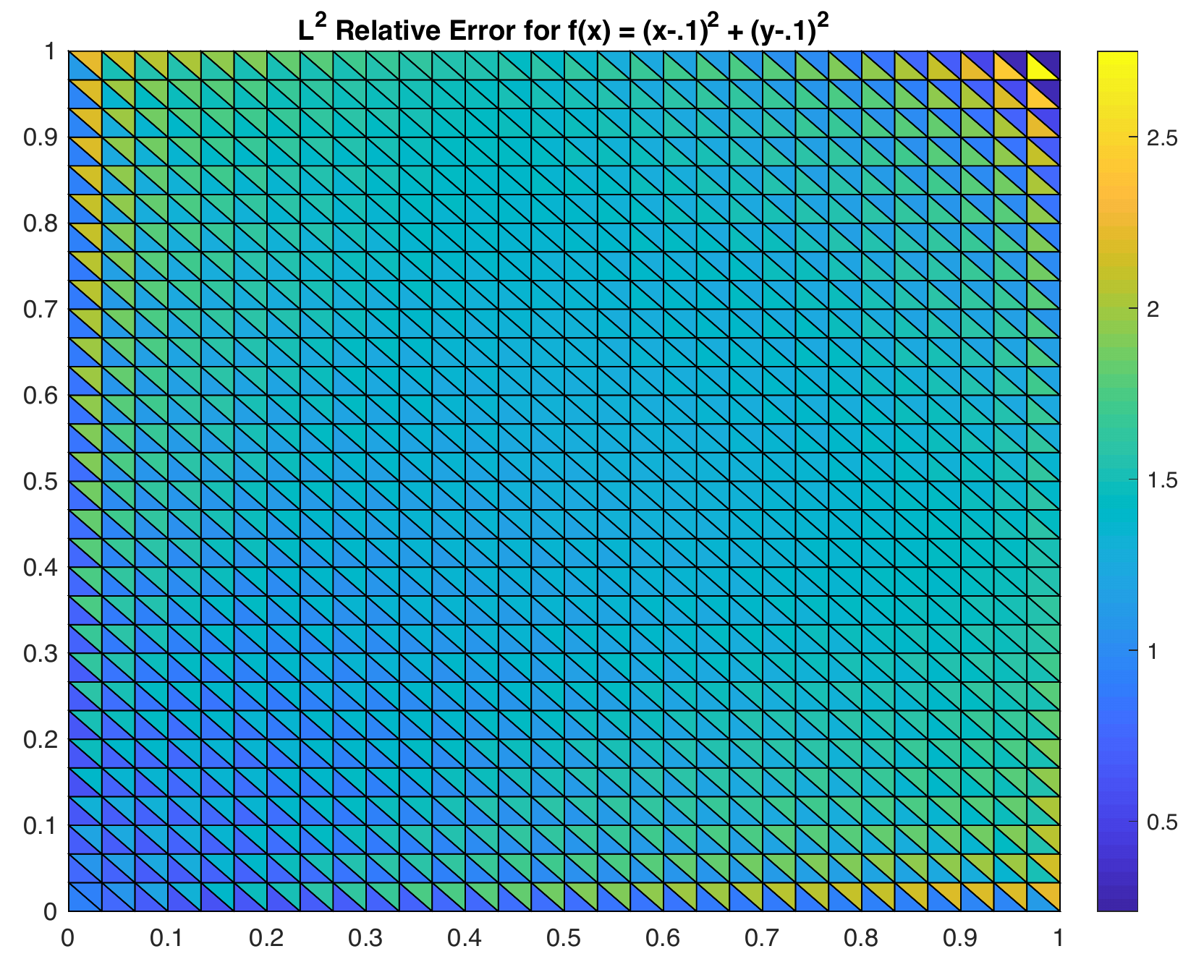}
        \caption{The triangles are colored to represent the relative $L^2$ error for $f(x) = (x-.1)^2 + (y-.1)^2$.  Errors here range from $.2369$ to $2.7504$, with an average error of $1.3561$.}\label{fig:RelErrorMinimum}
\end{figure}
\end{example}

\section{Discussion and conclusions}
\label{sec:disc}
In this paper we gave a new definition of the fractional discrete exterior derivative that satisfied three key properties.  Our definition was based on the fractional Caputo derivative, and we tested on some $1$- and $2$-d examples.  Our definition requires an interpretation of distance between two arbitrary $p$-simplices on the mesh.

There are some limitations in the fractional discrete exterior derivative definition.  First, it complicates the process of taking a discrete exterior derivative, which is normally represented by a single matrix with a nice structure.  While this does mean that the process of taking a fractional discrete exterior derivative is more involved, this is to be expected in comparison to a normal discrete exterior derivative.  Specifically, the non-local nature of the fractional derivative means that we have to dedicate time in the computation to determining distances between simplices that are far apart.   The other limitation is that our definition does not act on critical points the same way that the fractional Caputo derivative does.  This means that near critical points, we end up with a much higher error term than we would like.  

For future work, we have many avenues to investigate.  First, we would like to study the way that the fractional discrete exterior derivative responds to tweaking the diagonal of the matrix $W$. Furthermore, we would like to be able to understand the behavior better at the critical points of a function.

There is also a connection that we have not yet discussed.  In the past couple of decades, there has been work on defining an extension of exterior calculus to a fractional exterior calculus setting, e.g.~\cite{Tarasov,CS2001}.  This includes a definition of a fractional exterior derivative, with the form 
\[d^s = \sum_{i = 1} ^n \frac{\partial^s}{\partial {x^i}^s} [dx^i]^s. \]
We intentionally chose not to make use of definition this when defining our operator.  Directly discretizing this fractional exterior derivative would not fit in the ``discrete first'' approach that DEC takes, and discretizing it in a way analogous to how the discrete exterior derivative is discretized would require a notion of a fractional boundary.  One possible strength of discretizing with this formula is the possibility of ensuring $d^s \circ d^s = 0$, as is shown in~\cite{CS2001}.  Still, fractional exterior calculus as a whole is a relatively new field, and it is not clear which properties are essential to preserve when considering their discretization.

Beyond just investigating more of the properties of the definition, we want to extend into a more DEC related setting.  The results on an interval and a rectangular region subdivided into triangles are promising, however one strength of DEC is its ability to work on manifold-like simplicial complexes.  With a definition of the fractional discrete exterior derivative, we should be able to extend current work using DEC to FDEs on simplicial complexes that represent more intricate geometric objects.

In the vein of solving FDEs, one of the prominent operators in FDEs is the fractional Laplacian, $(-\Delta)^s$.  The Laplacian operator has a nice definition in DEC where $\Delta = * d * d$.  Hence studying the fractional Laplacian in the DEC setting may reduce to studying $ \star \mathbb{D}_{p-1}^s \star \mathbb{D}_{0}^s$, with $\star$ being some sort of fractional discrete Hodge star operator. A good overview of the work that has been done on the fractional Laplacian was written by Lischke et. al in \cite{Lischke2018}. 

As discussed in the introduction, there was also a shape analysis application by Paquet and Viktor in \cite{PV2013}  using the fractional de Rham operator.  In their approach, they take a fractional power by using an eigenvalue decomposition, and raising the eigenvalues to the fractional power.  In the fractional Laplacian literature, this is referred to as the spectral method for computing the fractional power of the operator.  This would be interesting to test against when we are ready to use our method to try discretizing the fractional Laplacian.  In general, we cannot use the spectral method for the fractional discrete exterior derivative because it is a non-square matrix.

We would also be interested in exploring what the other definitions of fractional derivatives yield when discretized through DEC.  For example, the Riemann-Liouville fractional derivative is given by

\begin{equation}
	^{RL}D_{[a,x]}^s f(x) = \dfrac{1}{\Gamma(n - s)} \dfrac{d^n}{dx^n}\int_a^x \dfrac{f(t)}{(x-t)^{n+s-1}} dt, 
\end{equation}

\noindent where $[a,b]$ is some domain of interest, $s \in \mathbb{R}^+$ and $n = \ceil s$.

The main difference between this definition and the Caputo definition is the order of the differentiation and integration.  Because the integration occurs before taking any derivatives, the class of functions we act on is slightly less restrictive.  We suspect that the Riemann-Liouville definition leads to a fractional discrete exterior derivative operator that has the form 

\begin{equation}
 	\frac{1}{\Gamma(1-s)} \mathbb{D}_p W\alpha
\end{equation}
\noindent where $W$ is a weighting matrix similar to before, but instead of acting on discrete $(p+1)$-forms, it acts on discrete $p$-forms.

Having a working version of the fractional discrete exterior derivative derived from the Riemann-Liouville derivative will be useful.  Many applications of FDEs require either the Caputo derivative or the Riemann-Liouville derivative, and they in general do not yield the same results.  Enabling the computation of both forms will allow for more tests that can be used to compare and evaluate fractional derivatives.  

We are also interested in using these operators to study fractional equivalents to problems that have been previously investigated with DEC.  For example, recent work has been done on spatial Darcy flow problems to numerically compare results against experimental data.  

\section*{Acknowledgements}
This material is based upon work supported by the U.S. Department of Energy, Office of Science, Office of Advanced Scientific Computing Research, under Award Number(s) DE-SC-0019039.

\section*{References}

\bibliography{Building_A_FDED}{}

\begin{thebibliography}{10}
\expandafter\ifx\csname url\endcsname\relax
  \def\url#1{\texttt{#1}}\fi
\expandafter\ifx\csname urlprefix\endcsname\relax\def\urlprefix{URL }\fi
\expandafter\ifx\csname href\endcsname\relax
  \def\href#1#2{#2} \def\path#1{#1}\fi

\bibitem{Hirani2003}
A.~N. Hirani, {Discrete Exterior Calculus}, Ph.D. thesis, Caltech (2003).

\bibitem{Dominitz2010}
A.~Dominitz, A.~Tannenbaum, Texture mapping via optimal mass transport, IEEE
  transactions on visualization and computer graphics 16~(3) (2010) 419--433.

\bibitem{Vaxman2016}
A.~Vaxman, M.~Campen, O.~Diamanti, D.~Panozzo, D.~Bommes, K.~Hildebrandt,
  M.~Ben-Chen, Directional field synthesis, design, and processing, in:
  Computer Graphics Forum, Vol.~35, Wiley Online Library, 2016, pp. 545--572.

\bibitem{Elcott2007}
S.~Elcott, Y.~Tong, E.~Kanso, P.~Schr{\"o}der, M.~Desbrun, Stable,
  circulation-preserving, simplicial fluids, ACM Transactions on Graphics (TOG)
  26~(1) (2007) 4.

\bibitem{mullen2011hot}
P.~Mullen, P.~Memari, F.~de~Goes, M.~Desbrun, {HOT}: {H}odge-optimized
  triangulations, in: ACM Transactions on Graphics (TOG), Vol.~30, ACM, 2011,
  p. 103.

\bibitem{ben2010discrete}
M.~Ben-Chen, A.~Butscher, J.~Solomon, L.~Guibas, On discrete killing vector
  fields and patterns on surfaces, in: Computer Graphics Forum, Vol.~29, Wiley
  Online Library, 2010, pp. 1701--1711.

\bibitem{zhang2010spectral}
H.~Zhang, O.~Van~Kaick, R.~Dyer, Spectral mesh processing, in: Computer
  graphics forum, Vol.~29, Wiley Online Library, 2010, pp. 1865--1894.

\bibitem{hormann2007mesh}
K.~Hormann, B.~L{\'e}vy, A.~Sheffer, Mesh parameterization: {T}heory and
  practice, presentation at SIGGRAPH (2007).

\bibitem{zaharescu2009surface}
A.~Zaharescu, E.~Boyer, K.~Varanasi, R.~Horaud, Surface feature detection and
  description with applications to mesh matching, in: Computer Vision and
  Pattern Recognition, 2009. CVPR 2009. IEEE Conference on, IEEE, 2009, pp.
  373--380.

\bibitem{leok2004foundations}
M.~Leok, Foundations of computational geometric mechanics, Ph.D. thesis,
  California Institute of Technology (2004).

\bibitem{mullen2011discrete}
P.~Mullen, A.~McKenzie, D.~Pavlov, L.~Durant, Y.~Tong, E.~Kanso, J.~E. Marsden,
  M.~Desbrun, Discrete {L}ie advection of differential forms, Foundations of
  Computational Mathematics 11~(2) (2011) 131--149.

\bibitem{kanso2007geometric}
E.~Kanso, M.~Arroyo, Y.~Tong, A.~Yavari, J.~G. Marsden, M.~Desbrun, On the
  geometric character of stress in continuum mechanics, Zeitschrift f{\"u}r
  angewandte Mathematik und Physik 58~(5) (2007) 843--856.

\bibitem{bridges2006numerical}
T.~J. Bridges, S.~Reich, Numerical methods for {H}amiltonian {PDE}s, Journal of
  Physics A: mathematical and general 39~(19) (2006) 5287.

\bibitem{Roop2006}
J.~P. Roop, {Computational aspects of FEM approximation of fractional advection
  dispersion equations on bounded domains in R2}, Journal of Computational and
  Applied Mathematics 193~(1) (2006) 243--268.
\newblock \href {http://dx.doi.org/10.1016/j.cam.2005.06.005}
  {\path{doi:10.1016/j.cam.2005.06.005}}.

\bibitem{Oktar2014}
O.~Ozgen, M.~Kallmann, L.~E. Ramirez, C.~F. Coimbra,
  \href{http://graphics.ucmerced.edu/papers/10-tog-fracdef.pdf}{Underwater
  cloth simulation with fractional derivatives}, ACM Transactions on Graphics
  (TOG) 29~(3) (2010) 1--9.
\newline\urlprefix\url{http://graphics.ucmerced.edu/papers/10-tog-fracdef.pdf}

\bibitem{Farquhar2018}
M.~E. Farquhar, T.~J. Moroney, Q.~Yang, I.~W. Turner, K.~Burrage,
  \href{http://arxiv.org/abs/1809.07936}{Computational modelling of cardiac
  ischaemia using a variable-order fractional laplacian}, arXiv:1809.07936
  (2018).
\newblock \href {http://arxiv.org/abs/1809.07936} {\path{arXiv:1809.07936}}.
\newline\urlprefix\url{http://arxiv.org/abs/1809.07936}

\bibitem{PV2013}
E.~Paquet, H.~L. Viktor, Isometrically invariant description of deformable
  objects based on the fractional heat equation, in: R.~Wilson, E.~Hancock,
  A.~Bors, W.~Smith (Eds.), Computer Analysis of Images and Patterns, Springer
  Berlin Heidelberg, Berlin, Heidelberg, 2013, pp. 135--143.

\bibitem{AMR2012}
R.~Abraham, J.~E. Marsden, T.~Ratiu, Manifolds, tensor analysis, and
  applications, Vol.~75, Springer Science \& Business Media, 2012.

\bibitem{Herrmann2010}
R.~Herrmann, {Fractional Calculus An Introduction for Physicists}, 2010.

\bibitem{Caputo1966}
M.~Caputo, Linear models of dissipation whose {Q} is almost frequency
  independent-{II}, Geophysical Journal International 13~(5) (1967) 529--539.

\bibitem{Tarasov}
V.~E. Tarasov, Fractional dynamics: applications of fractional calculus to
  dynamics of particles, fields and media, Springer Science \& Business Media,
  2011.

\bibitem{Mariya}
M.~K. Ishteva, {Properties and Applications of the Caputo Fractional Operator}
  (2005).

\bibitem{Garrappa2015}
R.~Garrappa, {Numerical evaluation of two and three parameter Mittag-Leffler
  functions}\href {http://arxiv.org/abs/1503.06569} {\path{arXiv:1503.06569}},
  \href {http://dx.doi.org/10.1137/140971191} {\path{doi:10.1137/140971191}}.

\bibitem{CS2001}
K.~Cottrill-Shepherd, M.~Naber, {Fractional differential forms}, Journal of
  Mathematical Physics 42~(5) (2001) 2203--2212.
\newblock \href {http://dx.doi.org/10.1063/1.1364688}
  {\path{doi:10.1063/1.1364688}}.

\bibitem{Lischke2018}
A.~Lischke, G.~Pang, M.~Gulian, F.~Song, C.~Glusa, X.~Zheng, Z.~Mao, W.~Cai,
  M.~M. Meerschaert, M.~Ainsworth, G.~E. Karniadakis,
  \href{http://arxiv.org/abs/1801.09767}{{What Is the Fractional Laplacian?}}
  (2018) 1--80\href {http://arxiv.org/abs/1801.09767}
  {\path{arXiv:1801.09767}}.
\newline\urlprefix\url{http://arxiv.org/abs/1801.09767}

\end{thebibliography}


\end{document}